\documentclass[11pt,reqno]{amsart}
\usepackage{fullpage}
\usepackage{dsfont}
\usepackage{amsmath}
\usepackage{amsfonts}
\usepackage{amssymb}
\usepackage{graphicx}
\usepackage{mathrsfs}
\usepackage{epsfig}
\usepackage{amsthm}
\usepackage{verbatim}
\usepackage{graphicx}
\usepackage{color}
\usepackage{url}

\usepackage[colorlinks,citecolor=black,linkcolor=black,
            bookmarksopen,
            bookmarksnumbered
           ]{hyperref}

\newcommand{\beq}{\begin{equation}}
\newcommand{\eeq}{\end{equation}}
\newcommand{\e}{\mathrm{e}}
\newcommand{\dd}{\,\mathrm{d}}

\newcommand{\taul}{\tau}
\newcommand{\vta}{u_{\tau}}
\newcommand{\uta}{u_{\tau}}

\newcommand{\unum}[1]{u^{#1}_\tau}
\newcommand{\un}{\unum{n}}

\newtheorem{theorem}{Theorem}[section]

\newtheorem{lemma}[theorem]{Lemma}
\newtheorem{cor}[theorem]{Corollary}

\newtheorem{rem}[theorem]{Remark}
\theoremstyle{definition}
\newtheorem{ex}[theorem]{Example}

\allowdisplaybreaks[3]

\author{Frédéric Rousset}
\address{Laboratoire de Mathématiques d’Orsay (UMR 8628), Université Paris-Saclay, CNRS, 91405 Orsay Cedex, France (F. Rousset)}
\email{frederic.rousset@universite-paris-saclay.fr}

\author{Katharina Schratz}
\address{LJLL (UMR 7598), Sorbonne Universit\'e, UPMC, 4 place Jussieu, 75005, Paris, France (K. Schratz)}
\email{katharina.schratz@sorbonne-universite.fr}
 
\begin{document}
\begin{abstract} A large toolbox of numerical schemes for dispersive equations has been established, based on different discretization techniques such as discretizing the variation-of-constants formula (e.g., exponential integrators) or splitting the full equation into a series of simpler subproblems (e.g., splitting methods). In many situations these classical schemes allow a precise and efficient approximation. This, however, drastically changes whenever non-smooth phenomena enter the scene such as for problems at low regularity and high oscillations. Classical schemes fail to capture the oscillatory nature of the solution, and this may lead to severe instabilities and loss of convergence. In this article we review a new class of resonance-based schemes. The key idea in the construction of the new schemes is to tackle and deeply embed the underlying nonlinear  structure of resonances into the numerical discretization. As in the continuous case, these terms are central to structure preservation and offer the new schemes strong properties at low regularity. 

\end{abstract}

\title[]{Resonances as a computational tool} 

\maketitle
\section{Introduction}
Nonlinear dispersive equations, e.g., the nonlinear Schr\"odinger, Korteweg--de Vries,   wave map equation, etc., have gained a lot of attention in the last decades. Their smooth solutions are nowadays well understood at the theoretical as well as computational level. While huge progress could also be made in their theoretical analysis for rough data, non-smooth solutions remain in large parts a mystery computationally. Unlike in parabolic problems we do not have strong smoothing effects. This leads to many interesting phenomena such as blow-up, growth of Sobolev norms, quantisation effects, turbulence, etc., which are, however extremely challenging to model numerically. As soon as roughness comes into play most classical numerical schemes break-down, and little is known on how to overcome this.

This survey article intends to review recent developments in so-called resonance based schemes which try to build a bridge between smooth and non-smooth numerics for  dispersive equations. A first attempt  of so-called \textcolor{black}{resonance-based schemes}  (see for instance \cite{BBS23,YYMS,BMS,BBS,BS,CRS22,CS22,HS16,ORS21,ORS22,OS18,RS,SWZ21}),   was  profoundly inspired by   theoretical analysis of  dispersive equations  at low regularity  (Bourgain \cite{Bour93a}, Tao~\cite{Tao06}) and rough path theory (Gubinelli~\cite{Gub11}) and provides  a powerful tool which in many situations allows for approximations in a much more general setting (i.e., for rougher data) than classical schemes  (e.g., Splitting methods),  see also the recent important works \cite{BLW,CLL23,LW21,LW22,WZ22KG,Mass,WZ22}.

The central aim of this survey article is to present in detail the main idea behind the novel technique on the basic test example of cubic, periodic Schr\"odinger equation (Section \ref{sec:idea}), provide various examples of dispersive equations and beyond (Section \ref{sec:general}), explain the difficulties in  establishing low-regularity error estimates (Section \ref{sec:error}), and outline open problems in this direction (Section \ref{sec:open}).

\section{The main idea}\label{sec:idea}
To explain the key idea behind {resonances as a computational tool} we first have to understand why  classical schemes (e.g.,  splitting methods and exponential integrators) in general fail to approximate  rough dynamics in dispersive equations. For this purpose let us consider as a simple model problem the one-dimensional, periodic cubic nonlinear Schr\"odinger (NLS) equation 
\begin{equation}\label{cubNLS}
i \partial_t u(t,x) = - \Delta u(t,x) + \vert u(t,x)\vert^2 u(t,x)\qquad (t,x)\in \mathbb{R}\times \mathbb{T}
\end{equation}
with rough initial data
\begin{equation}\label{iniSig}
u(0,x) = u_0(x) \in H^\sigma(\mathbb{T}),
\end{equation}
where we want to choose $\sigma >0$ as small as possible (we will see later in Section \ref{sec:error} how far we can actually push down the Sobolev index $\sigma$ to obtain convergence estimates in $L^2(\mathbb{T})$). Here, $H^\sigma(\mathbb{T})$ denotes the classical Sobolev space on   torus $\mathbb{T}$ with regularity $\sigma\geq 0$.\vskip0.2cm

\subsection{Splitting methods} One of the most famous  numerical methods to approximate the time dynamics of linear and nonlinear Schr\"odinger equations  such as \eqref{cubNLS} are splitting methods (see, e.g., \cite{Long1,Faou12,Birk1,Birk2,HLW,HLRS10,HUV,HLRS10,Lubich08,McLacQ02,Shi}). The main idea  lies in splitting the full equation into a series of simpler subproblems. One then solves these subproblems (either exactly or with a numerical scheme) and composes the sub-flows to obtain an approximation to the flow of the original -- full -- equation.  Looking at the Schr\"odinger equation \eqref{cubNLS} we face two main challenges numerically: The differential operator $-\Delta$ and the non-linearity $\vert u\vert^2u$. Instead of solving the full problem \eqref{cubNLS} the idea of splitting for NLS is to consider its kinetic (T) and nonlinear part (V) separetely, i.e., 
\begin{align}
\text{ (T) } \quad i \partial_t u_T = -\Delta u_T \quad \text{and} \quad\text{ (V) } \quad i \partial_t u_V = \vert u_V\vert^2 u_V.
\end{align}
The  main advantage lies in the fact that both subproblems can be solved exactly in time: (T) in Fourier space (in case of a spatial discretisation, with a Fourier pseudo spectral method) and (V) as the modulus $\vert u_V\vert^2$ is conserved, see for instance \cite{Faou12,Lubich08} for details. A simple composition of the kinetic and nonlinear sub-flow   leads at first-order  to the so-called Lie splitting method which at time $t_n = n \tau$ ($\tau$ denoting the time step size) takes the form
\begin{equation}\label{Lie}
u^{n+1} = \underbrace{e^{ i \tau \Delta}}_{\text{(T)}} \underbrace{e^{-i \tau \vert u^n\vert^2} u^n}_{\text{(V)}}.
\end{equation}
The local error of Lie splitting \eqref{Lie} is driven by the Lie commutator $[T,V](u)$  which reads
\begin{equation}\label{commLie}
\frac{1}{2} [T,V](u) =  \left(\nabla u \cdot \nabla \overline{u}\right) u + \left(\nabla u \overline{u} \right)\cdot \nabla u + \left(u \nabla \overline{u}\right)\cdot \nabla u + \left(u \overline{\Delta u} \right)u,
\end{equation}
see \cite[Section 4.2]{Lubich08}. Due to the appearance of $ \overline{\Delta u}$ in the local error, the boundedness of at least two additional derivatives of the exact solution is required.  Higher-order splitting methods require more   regularity as they  introduce a local error with nested commutators. In case of second-order Strang splitting the error is for instance driven by the double commutator $[T,[T,V]]$  involving the term $ \overline{\Delta \Delta u}$. This requires the boundedness of four additional derivatives. For a precise convergence analysis we refer to \cite{Lubich08}.


\subsection{Exponential integrators} 
Another well-known method to solve Schr\"odinger-type equations are exponential integrators  (see, e.g., \cite{CCO08,HLW,HoLu99,HochOst10} and the references therein). The main idea   lies in discretising   Duhamel's formula which for cubic NLS \eqref{cubNLS}  takes the form
\begin{equation}\label{Duh}
u(t) = e^{it\Delta} u(0) - i e^{i t \Delta} \int_0^t e^{-i s\Delta}\left( \vert u(s)\vert^2 u(s)\right) \mathrm{d}s.
\end{equation}
At  time $t_{n+1}=t_n + \tau$ we find (considering Duhamels formula on the time interval $[0,\tau]$  with  initial value $u(t_n)$) that
 \begin{equation}\label{DuhN}
u(t_{n+1}) = e^{i \tau \Delta} u(t_n) - i e^{i \tau \Delta} \int_0^\tau e^{-i s\Delta}\left(\vert u(t_n+s)\vert^2 u(t_n+s) \right)\mathrm{d}s.
\end{equation}
Exponential integrator schemes are based on  Taylor series expansion of the solution  within the integral. At first order they build on the first-order Taylor series expansion 
\begin{equation}\label{taylor1exp}
u(t_n+s) = u(t_n) + \mathcal{O}( s u').
\end{equation}
Plugging the approximation \eqref{taylor1exp} into Duhamel's formula \eqref{DuhN} leads to the first-order exponential integrator scheme
\begin{equation}\label{expInt}
u^{n+1} = e^{i \tau \Delta} u^n - i \tau \varphi_1(i \tau \Delta) \vert u^n\vert^2 u^n \quad \text{with}\quad\varphi_1(z) = \frac{e^{z}-1}{z}.
\end{equation}
From the Taylor series expansion \eqref{taylor1exp} we easily see that the local error of the first-order exponential integrator method \eqref{expInt} is driven by the time derivative $u'$, where 
\begin{equation}\label{errExp}
\mathcal{O}(u') = \mathcal{O}(\Delta u)
\end{equation}
in sense of derivatives. Hence, as for Lie splitting first-order convergence requires the boundedness of at least two additional derivatives. Higher-order exponential integrators are based on higher-order Taylor series expansion of the solution within the integral, i.e., 
at second-order one takes
\begin{equation}\label{taylor2exp}
u(t_n+s) = u(t_n) +s u'(t_n) + \mathcal{O}(s^2 u'')
\end{equation}
and replaces the time derivative $u'(t_n)$ by the equation itself $u'(t_n) = i \Delta u(t_n) - i \vert u(t_n)\vert^2 u(t_n)$.   Due to the local error scaling
$$
 \mathcal{O}(u'')= \mathcal{O}(\Delta\Delta u)
$$
we see that as for second-order Strang splitting, the second-order exponential integrator method requires the boundedness of at least four additional derivatives.

\subsection{Structure of the solution and classical methods}
From  local error structure \eqref{commLie} and~\eqref{errExp} we see that both the Lie splitting and exponential integrator method require smooth solutions: At least the boundedness of two additional derivatives is necessary for their first-order convergence. Their second-order counterparts require even the boundedness of four additional derivatives, and in general we have a local error scaling at order $\nu$ of type
\begin{equation}\label{locNu}
\mathcal{O}(\tau^{\nu+1} (-\Delta)^\nu)
\end{equation}
such that for instance for a method of order four we would  require the boundedness of 16 additional derivatives. A natural question therefore arises: {\em Can we construct numerical schemes which allow  convergence for rougher data than classical methods, i.e., under lower regularity assumptions than splitting or exponential integrators?}\smallskip

To answer this question we first have to understand the underlying structure of the solution in a better way. For this purpose we turn back to Duhamel's formula \eqref{Duh}. However, instead of approximating its solution $u(s)$ by a classical Taylor series expansion, we look at iterations of   Duhamel's formula: Using that
$$
u(s) = e^{i s \Delta}u(0) + \int_0^s \ldots \mathrm{d}s_1
$$
we find
\begin{equation}\label{DuhIt}
u(t) = e^{it\Delta} u(0) - i e^{i t \Delta} \int_0^t e^{-i s\Delta}\left( \left\vert  e^{is\Delta} u(0) \right\vert^2 e^{is\Delta} u(0)\right)\mathrm{d}s +  \int_0^t  \int_0^s \ldots \mathrm{d}s_1\mathrm{d}s.
\end{equation}
If we want to get a rough idea of the dynamics of the solution, we can at first forget about the higher order iterations, and neglect the double integral $ \int_0^t  \int_0^s\mathrm{d}s_1\mathrm{d}s$ in \eqref{DuhIt}. We then see that the underlying structure of the solution $u(t)$ is driven by the nonlinear frequency interaction of $-\Delta$ and $\pm \Delta$ with leading oscillations
\begin{equation}\label{osc}
\mathrm{Osc}(s,\Delta, u(0)) = e^{-i s\Delta}\left( \left\vert  e^{is\Delta} u(0) \right\vert^2 e^{is\Delta} u(0)\right) 
.
\end{equation}
Numerical schemes stay close to the structure of the solution -- even at low regularity --  if they   resolve  --  even for rough data -- the leading oscillations~\eqref{osc}.

A closer look, however, shows that splitting methods and exponential integrators in general neglect the nonlinear frequency interactions in \eqref{osc}: Lie splitting \eqref{Lie}  is based on the frequency approximation
 \begin{equation}\label{Lie:freq}
\mathrm{Osc}(s,\Delta, u(0))  \approx \vert u(0)\vert^2 u(0),
\end{equation}
while exponential integrator methods  swallow all frequencies within the nonlinearity based on  \begin{equation}\label{Exp:freq}
\mathrm{Osc}(s,\Delta, u(0))  \approx 
e^{-i s\Delta}\left( \left\vert    u(0) \right\vert^2 u(0)\right).
\end{equation}
In case of smooth solutions, for which $\Delta u$   is well defined in the space of interest (see also error structure \eqref{commLie} and \eqref{errExp}),   linearisation of frequencies such as \eqref{Lie:freq} and \eqref{Exp:freq}   in general lead to good approximations of the exact solution \eqref{DuhIt}. This can be seen by a simple Taylor series expansion of the oscillations
\begin{equation}\label{taylor:class}
e^{\pm i s\Delta} u = u + \mathcal{O}(s \Delta u).
\end{equation}
Expansion \eqref{taylor:class} introduces a small remainder of order $s$ as long as $u$ is sufficiently smooth, i.e., $\Delta u$ is bounded. For rough solutions, for which 
$$
\Delta u
$$
becomes unbounded,   approximations such as \eqref{Lie:freq} and \eqref{Exp:freq}  in general, however, break down as the linearisation of frequencies \eqref{taylor:class} is no longer valid. The latter is not only a theoretical artefact stemming from error analysis, but also drastically observed in numerical experiments.\smallskip

The main aim of {resonances as a computational tool} is to  overcome this by stepping away from linearised frequency interactions (such as   \eqref{Lie:freq} and \eqref{Exp:freq}) towards new schemes which deeply embed the nonlinear frequency interactions  of the solution  (or at least their dominant parts) into the numerical discretisation. In general this allows for much rougher data than classical schemes.  In the following we will explain this idea in detail on the concrete example of cubic, periodic NLS~\eqref{cubNLS}. We will then give a more general overview of the new ansatz for a broad class of equations in Section~\ref{sec:general}.

\subsection{Resonances as a computational tool: The main idea}\label{subsec:rc1}
In case of cubic, periodic  NLS~\eqref{cubNLS} the underlying structure of the solution $u$ is driven be the leading oscillations  \eqref{osc}. We aim to find a good approximation to them, even for rough data. For this purpose we turn to Fourier analysis as this will allow us to analyse exactly  the nonlinear frequency interactions of $\Delta$ and $\pm \Delta$ in \eqref{osc}: With $\hat{u}_k$ denoting the $k-$th Fourier coefficient in space, i.e.,  
$$u(t,x) = \sum_{k\in \mathbb{Z}} e^{i k x} \hat{u}_k(t)$$  the central oscillations \eqref{osc} of cubic NLS  take the form
\begin{equation}\label{oscNLSf}
\mathrm{Osc}(s,\Delta, u(0)) = \sum_{k = k_1 -k_2+k_3}e^{i k x} \hat{u}_{k_1}(0)\hat{u}_{k_2}(0)\overline{\hat{u}}_{k_3}(0)
 e^{i s k^2 } e^{-isk_1^2}  e^{isk_2^2} e^{-isk_3^2} .
 \end{equation}
Hence, we find that cubic NLS \eqref{cubNLS} is driven by the underlying   resonance structure
 \begin{equation}\label{freq:NLS}
\mathrm{R}(k_1,k_2,k_3) =  k^2 -k_1^2 + k_2^2 - k_3^2 
 \end{equation}
 which determines the nonlinear frequency interactions in \eqref{oscNLSf}. \smallskip
 
The problem of  classical methods is that they in general linearise the nonlinear   resonance structure \eqref{freq:NLS} and treat -- from a frequency point of view -- the nonlinear PDE \eqref{cubNLS} as if it was a linear problem. More precisely,  
Lie splitting \eqref{Lie} approximates $\mathrm{R}$ by zero
$$
\mathrm{R}(k_1,k_2,k_3)  \approx 0 \quad \text{for all}\quad k_1,k_2, k_3 \in \mathbb{Z}
$$
 (cf. \eqref{Lie:freq}). Exponential integrator methods \eqref{expInt}, on the other hand, swallow the frequencies $-k_1^2, k_2^2$ and $ k_3^2$ within the nonlinearity and approximate the resonance structure $\mathrm{R}$ as follows 
$$
\mathrm{R}(k_1,k_2,k_3)  \approx k^2\quad \text{for all}\quad k_1,k_2, k_3 \in \mathbb{Z}
$$
(cf. \eqref{Exp:freq}). The central idea in resonances as a computational tool is to step away from linearised frequency approximations towards an improved nonlinear approach: Instead of linearising the resonance structure $\mathrm{R}(k_1,k_2,k_3)$ we filter out its dominant part, and solve the dominant part exactly, while  only approximating the lower-order parts. The central question is: {\em What is actually the dominant part in the nonlinear resonance structure \eqref{freq:NLS}? And is there a unique way to define it?}

It turns out that there are many ways to define the dominant part with a lot of questions remaining open so far (see also Section \ref{sec:open}). Here we focus on classical resonance-based methods, which see dominance in terms of derivatives:

 \begin{rem}\label{rem:domNLS}
 If we take a closer look at the resonance structure  \eqref{freq:NLS} we see that $\mathrm{R}$ can be expressed as follows:
 \begin{equation}\label{Res:NLS}
\mathrm{R}(k_1,k_2,k_3) =  k^2 -k_1^2 + k_2^2 - k_3^2  = 2 k_2^2 - 2 (k_1+k_3) k_2 + 2 k_1 k_3.
 \end{equation} 
We observe that the mixed terms in \eqref{Res:NLS}
$$
- 2 k_1k_2, \quad -2 k_3 k_2, \quad 2 k_1 k_3
$$
correspond to first-order derivatives, while the quadratic term
$$
2k_2^2
$$
corresponds to a second-order derivative. This can  easily be seen in Fourier space: Let us take two smooth functions $v$ and $w$
$$
v(x) = \sum_{\ell \in \mathbb{Z}} \hat{v}_\ell e^{i \ell x}, \quad 
w(x) = \sum_{m \in \mathbb{Z}} \hat{w}_m e^{i m x}.
$$
Then we have that 
$$
 \partial_x v \cdot \partial_x w = - \sum_{\ell, m \in \mathbb{Z}} \ell \cdot m \,
 \hat{v}_\ell\hat{w}_m e^{i (\ell+m) x} \quad \text{while } \quad w \partial_x^2 v = - \sum_{\ell, m \in \mathbb{Z}} \ell^2\,
 \hat{v}_\ell\hat{w}_m e^{i (\ell+m) x}.
$$
\end{rem}
Thanks to Remark \ref{rem:domNLS} we see that the dominant part, with the highest order of derivative, in the nonlinear frequency interactions \eqref{Res:NLS} is the quadratic term
$$
2 k_2^2.
$$
Hence, the idea is to treat the second-order term $2 k_2^2$ exactly within the numerical discretisation, and only approximate the lower-order mixed terms $$- 2 k_1k_2, \quad-2 k_3 k_2,\quad 2 k_1 k_3.$$ This can be achieved by Taylor series expansion
\begin{equation}\label{Tt}
 e^{i s k^2 } e^{-isk_1^2}  e^{isk_2^2} e^{-isk_3^2}  = 
 e^{i s \left( 2 k_2^2 - 2 (k_1+k_3) k_2 + 2 k_1 k_3\right)} = 
  e^{ 2 i s   k_2^2 } +\mathcal{O}\big(s 
 \left[2 (k_1+k_3) k_2 + 2 k_1 k_3\right] 
 \big).
\end{equation}
The main advantage of the resonance-based approximation \eqref{Tt} lies in the fact that it introduces a local error involving only first-order instead of  second-order derivatives (cf. \eqref{taylor:class}), building the basis of the first-order resonance-based approximation
\begin{equation}\label{oscResapp}
\begin{aligned}
\mathrm{Osc}(s,\Delta, u(0))  & = u^2(0) e^{-2is\Delta} \overline u(0)+ \mathcal{O}\left(s  \left(\partial_x u\right) \partial_x (u \overline u )\right).
 \end{aligned}
 \end{equation}
Plugging the  approximation \eqref{oscResapp} into the iteration of Duhamel's formula \eqref{DuhIt} yields that
\begin{equation}\label{app:res}
\begin{aligned}
u(t) & = e^{it\Delta} u(0) - i e^{i t \Delta} \int_0^t u^2(0) e^{-2is\Delta} \overline u(0) \mathrm{d}s +   \mathcal{O}\left(t^2  \left(\partial_x u\right) \partial_x (u \overline u )\right)
+ \mathcal{O}\left(t^2  u^5\right)\\
& = e^{it\Delta} u(0) - i t e^{i t \Delta} \left( u^2(0) \varphi_1\left(-2it\Delta\right) \overline u(0) \right)+   \mathcal{O}\left(t^2  \left(\partial_x u\right) \partial_x (u \overline u )\right)
+ \mathcal{O}\left(t^2  u^5\right).
\end{aligned}
\end{equation}
The resonance-based approximation \eqref{app:res} allows for rougher data than classical approximations (e.g., Lie splitting \eqref{Lie} and exponential integrator \eqref{expInt}) which are based on the Taylor-series expansion \eqref{taylor:class}: Instead of the full Laplacian, the improved remainder term \eqref{app:res} only involves first-order derivatives.

\begin{rem}
The resonance-based approximation of Duhamel's formula \eqref{app:res} leads to the following numerical scheme
$$
u^{n+1} = e^{i\tau \Delta} u^n - i \tau e^{i \tau \Delta} \left( (u^n)^2 \varphi_1\left(-2i\tau\Delta\right) \overline {u^n} \right)
$$
which was originally introduced in \cite{OS18}. This schemes is of first order, due to an approximation error at order $\tau^2$ (cf. \eqref{app:res}). For the error analysis  we refer to Section \ref{sec:error}. 
\end{rem}
\begin{rem}\label{rem:fourier}
One might now ask what happens if we want to construct numerical schemes which do not require any additional smoothness in the solution at the time discrete level (for the spatial discretisation some regularity is in general always needed). In theory this would be possible by resolving all nonlinear frequency interactions  \eqref{Res:NLS} exactly within the numerical discretisation. For certain particular one-dimensional equations this is indeed possible (such as the KdV equation \cite{HS16,LW22,WZ22} and periodic, one-dimensional cubic NLS \cite{Mass}). In general, however,  practical implementation of fully resonance-based schemes in higher dimensions would lead to huge computational costs as all computations would need to be carried out in Fourier space. The aim of resonances as a computational tool  lies in constructing a new class of schemes which allow for rougher data than classical methods, but at similar computational costs. Thus, we want to carry out differentation in Fourier space, and function multiplication in physical space. This allows for fast computations with the aid of the discrete Fourier transform. 
\end{rem}

\subsection{Resonances as a computational tool: What about higher-order?}\label{subsec:rc}
A next natural question is whether we can use this idea to achieve higher-order schemes for rougher data than classical methods require, i.e., reduce the regularity assumptions in the classical error scaling \eqref{locNu}. Two steps are essential: 
\begin{itemize}
\item[(I)] The iteration of Duhamel's formula \eqref{Duh} up to higher order, and
\item[(II)] the approximation of leading oscillations  \eqref{osc} up to higher order.
\end{itemize}
We will address these points in the following two subsections separately. Interestingly, it turns out that (II) is in fact much harder to accomplish than (I).
\subsubsection{Higher-order iteration of Duhamel's formula}
Let us consider the second-order iteration of Duhamel's formula  \eqref{Duh}. To obtain our resonance-based numerical approximation at second order  we take the  first-order resonance-based approximation~\eqref{app:res} at time $t=s$
\begin{equation*}\begin{aligned}
u(s) 
& = e^{is\Delta} u(0) - i s e^{i s \Delta} \left( u^2(0) \varphi_1\left(-2is\Delta\right) \overline u(0) \right)+   \mathcal{O}\left(s^2  \left(\partial_x u\right) \partial_x (u \overline u )\right)
+ \mathcal{O}\left(s^2  u^5\right)
\end{aligned}
\end{equation*}
and plug it   into Duhamel's formula \eqref{Duh}. This leads to the following expansion of the solution
\begin{equation}\label{Duh2}
\begin{aligned}
u(t) &= e^{it\Delta} u(0) - i e^{i t \Delta} \int_0^t e^{-i s\Delta}\left( \vert  e^{is\Delta} u(0) \vert^2  e^{is\Delta} u(0) \right) \mathrm{d}s\\
& - (-2 i) i e^{i t \Delta} \int_0^t  s e^{-i s\Delta}\Big[\vert  e^{is\Delta} u(0) \vert^2   e^{i s \Delta} \left( u^2(0) \varphi_1\left(-2is\Delta\right) \overline u(0) \right)  \Big] \mathrm{d}s\\
& -(i)  i e^{i t \Delta} \int_0^t  s e^{-i s\Delta}\Big[ (  e^{is\Delta} u(0) )^2  e^{-i s \Delta} \left( \overline u^2(0)\overline{\varphi_1\left(-2is\Delta\right)} u(0) \right) \Big] \mathrm{d}s\\
& +    \mathcal{O}\left(t^3  \partial_x u \right),
\end{aligned}
\end{equation}
where $ \mathcal{O}\left(t^3  \partial_x u \right)$ denotes a remainder of polynomials in $u, \overline u, \partial_x u$ and $\partial_x \overline u$.

\begin{rem}\label{rem:Reg2} If we want to reach second-order accuracy -- so a local error at order $\tau^3$ -- we have to sacrifice in terms of regularity. While classical second-order schemes such as Strang splitting  or second-order exponential integrators introduce a local  error at order
$$
\mathcal{O}(\tau^3 \Delta^2 u)
$$
(cf. \eqref{locNu})  a resonance-based approach allows us to halve these regularity assumptions with a local error structure of type
 \begin{equation}\label{err2}
\tau^3 \Delta u.
\end{equation}
It remains an open question whether these regularity assumptions can be reduced even further, and second-order convergence can be reached under the same regularity assumptions as first-order.
\end{rem}
Error structure \eqref{err2} allows for the loss of two derivatives in the remainder. Hence, we can expand~\eqref{Duh2} even  further: Using the classical Taylor series expansion $e^{\pm i s \Delta} = 1 + \mathcal{O}(s \Delta)$ in the second and third integral in \eqref{Duh2} we obtain
\begin{equation}\label{Duh3}
\begin{aligned}
u(t) &= e^{it\Delta} u(0) - i e^{i t \Delta} \int_0^t e^{-i s\Delta}\left( \vert  e^{is\Delta} u(0) \vert^2  e^{is\Delta} u(0) \right) \mathrm{d}s\\
& -2  e^{i t \Delta} \int_0^t  s  \vert   u(0) \vert^2    \left( u^2(0) \varphi_1\left(-2is\Delta\right) \overline u(0) \right)  \mathrm{d}s\\
& +  e^{i t \Delta} \int_0^t  s \,     u^2(0)    \left( \overline u^2(0)\overline{\varphi_1\left(-2is\Delta\right)} u(0) \right) \mathrm{d}s\\
& +    \mathcal{O}\left(t^3  \Delta u \right),
\end{aligned}
\end{equation}
where $ \mathcal{O}\left(t^3  \Delta u \right)$ denotes a remainder in polynomials at most of order $\Delta u$. Next we calculate that
\begin{align*}
& \int_0^t s \varphi_1\left(-2i s \Delta  \right)   \overline u(0) \mathrm{d}s =
 \frac1{-2i \Delta} \int_0^t \left(e^{-2is \Delta} - 1\right) \overline u(0)\mathrm{d}s
 =  \frac{t^2}{2}   \overline u(0) +  \mathcal{O}\left( t^3  \overline u(0)  \right) \\
 & \int_0^t s \overline{\varphi_1\left(-2i s \Delta  \right)} u(0)\mathrm{d}s
 =
 \frac1{2i \Delta} \int_0^t \left(e^{2is \Delta} - 1\right)u(0)\mathrm{d}s
 =  \frac{t^2}{2}  u(0) +  \mathcal{O}\left( t^3   u(0)  \right).
\end{align*}
Plugging this into iteration of Duhamel's formula \eqref{Duh3} we obtain
\begin{equation}\label{Duh4}
\begin{aligned}
u(t) &= e^{it\Delta} u(0) - i e^{i t \Delta} \int_0^t e^{-i s\Delta}\left( \vert  e^{is\Delta} u(0) \vert^2  e^{is\Delta} u(0) \right) \mathrm{d}s  -\frac{t^2}{2} \vert u(0)\vert^4 u(0)  +    \mathcal{O}\left(t^3  \Delta u \right),
\end{aligned}
\end{equation}
where $ \mathcal{O}\left(t^3  \Delta u \right)$ denotes a remainder in polynomials at most of order $\Delta u$ (which goes in line with Remark \ref{rem:Reg2}).
\begin{rem}
Thanks to the expansion \eqref{Duh4} we see that the main difficulty in achieving higher order resonance-based schemes lies in finding a suitable approximation of the leading oscillatory integral (cf. \eqref{osc})
\begin{equation}\label{intOsciL}
\int_0^t \mathrm{Osc}(s,\Delta, u(0))  \mathrm{d}s = \int_0^t e^{-i s\Delta}\left( \vert  e^{is\Delta} u(0) \vert^2  e^{is\Delta} u(0) \right) \mathrm{d}s .
\end{equation}
\end{rem}
\subsubsection{Higher-order approximation of  the oscillatory integral \eqref{intOsciL}}
In Fourier space  it holds that (see also~\eqref{oscNLSf} and \eqref{Res:NLS})
\begin{equation}\begin{aligned}\label{osciF2}
\mathrm{Osc}(s,\Delta, u(0)) 
 & = \sum_{k = k_1 -k_2+k_3}e^{i k x}\hat{u}_{k_1}(0)\overline{\hat{u}}_{k_2}(0)\hat{u}_{k_3}(0)
 e^{i s \text{R}(k_1,k_2,k_3)}
 \end{aligned}
 \end{equation}
 with
 $$
  \text{R}(k_1,k_2,k_3) = 2 k_2^2 - 2 (k_1+k_3) k_2 + 2 k_1 k_3.$$
Thanks to Remark \ref{rem:domNLS} we identify $2 k_2^2$ as the dominant part in the nonlinear resonance structure~$ \text{R}$ and can thus approximate the mixed terms $- 2 k_1 k_2,  - 2 k_3 k_2$ and $2 k_1 k_3$. At first-orderwe could simply carry out a first-order Taylor series expansion in the lower-order   terms cf. \eqref{Tt}. To achieve a higher-order resonance-based approximation the natural   idea is thus to apply a higher-order Taylor-series expansion in the lower order terms, i.e.,
\begin{equation}\label{TaylorF2}
e^{i s\text{R}(k_1,k_2,k_3)} = 
e^{2i s   k_2^2} - 2i s  \left(- 2 (k_1+k_3) k_2 + 2 k_1 k_3\right) + 
 \mathcal{O}\left( s^2 \left(- 2 (k_1+k_3) k_2 + 2 k_1 k_3\right)^2\right)
\end{equation}
and to neglect the terms of order $s^2$ involving at most second-order derivatives (which goes in line with Remark \ref{rem:Reg2}). The term
$$
2i s \left(- 2 (k_1+k_3) k_2 + 2 k_1 k_3\right) 
$$
is, however, unbounded, and including it in our numerical discretisation would  lead to loss of stability.

 To overcome this,  we choose a different route: We define the lower-order frequency interactions as follows
$$
 \psi_{\text{low}}(s, k_1,k_2,k_3) = e^{i s \left( - 2 (k_1+k_3) k_2 + 2 k_1 k_3\right)}
$$
which allows us to explicitly single out the dominant oscillations $ e^{2i s   k_2^2 } $ in \eqref{osciF2} thanks to the representation
\begin{equation*}\begin{aligned}
\mathrm{Osc}(s,\Delta, u(0)) 
 & = \sum_{k = k_1 -k_2+k_3}e^{i k x}\hat{u}_{k_1}(0)\overline{\hat{u}}_{k_2}(0)\hat{u}_{k_3}(0)\,
 e^{2i s   k_2^2 } \cdot  \psi_{\text{low}}(s, k_1,k_2,k_3) .
 \end{aligned}
 \end{equation*}
 Next we use a {\em stabilised} Taylor-series expansion to approximate the lower-order parts $\psi_{\text{low}}$: For $0 \leq s \leq t$ we have
 \begin{equation}\label{stab:TaylorF2}
  e^{2i s   k_2^2 } \psi_{\text{low}}(s, k_1,k_2,k_3) =   e^{2i s   k_2^2 } \left(1 + s \frac{  \psi_{\text{low}}(t, k_1,k_2,k_3)  -  \psi_{\text{low}}(0, k_1,k_2,k_3) }{t} + \mathcal{O}(t^2 \psi_{\text{low}}'') \right),
  \end{equation}
where $ \psi_{\text{low}}''$ involves at most second order derivatives. Plugging the stabilised Taylor-series expansion~\eqref{stab:TaylorF2} into \eqref{intOsciL} we obtain
\begin{equation*}\begin{aligned}
 \int_0^t & \mathrm{Osc}(s,\Delta, u(0))\mathrm{d}s  \\ &= \int_0^t \sum_{k = k_1 -k_2+k_3}e^{i k x} \hat{u}_{k_1}(0)\overline{\hat{u}}_{k_2}(0)\hat{u}_{k_3}(0)
e^{2i s   k_2^2 } \left(1 + s \frac{ \psi_{\text{low}}(t, k_1,k_2,k_3)  -  \psi_{\text{low}}(0, k_1,k_2,k_3) }{t} \right)\mathrm{d}s
 \\ &+ \mathcal{O}(t^3 \Delta u)\\
 &=  \sum_{k = k_1 -k_2+k_3}e^{i k x} \hat{u}_{k_1}(0)\overline{\hat{u}}_{k_2}(0)\hat{u}_{k_3}(0)
 \left(t\varphi_1{(2i t   k_2^2 )} +t^2 \varphi_2{(2i t   k_2^2 )} \frac{ \psi_{\text{low}}(t, k_1,k_2,k_3)  -   \psi_{\text{low}}(0, k_1,k_2,k_3) }{t} \right) 
  \\& + \mathcal{O}(t^3 \Delta u).
 \end{aligned}
 \end{equation*}
Now we use that
$$
\sum_{k_2}e^{- i k_2 x}  \overline{\hat{u}}_{k_2}(0)  \varphi_\sigma{(2i t   k_2^2 )}
=  \varphi_\sigma\left(-2i t\Delta \right) \overline{u(0)}
$$
as well as the relation
\begin{align*}
 \sum_{k = k_1 -k_2+k_3}e^{i k x} &\hat{u}_{k_1}(0)\overline{\hat{u}}_{k_2}(0)\hat{u}_{k_3}(0)
\psi(t, k_1,k_2,k_3)
\\&  = \sum_{k = k_1 -k_2+k_3}e^{i k x} \hat{u}_{k_1}(0)\overline{\hat{u}}_{k_2}(0)\hat{u}_{k_3}(0) e^{i t \left( 2 k_2^2 - 2 (k_1+k_3) k_2 + 2 k_1 k_3\right)} e^{-2 it k_2^2}\\
& = \sum_{k = k_1 -k_2+k_3}e^{i k x}\hat{u}_{k_1}(0)\overline{\hat{u}}_{k_2}(0)\hat{u}_{k_3}(0) e^{it k^2 } e^{-itk_1^2}  e^{-itk_2^2} e^{-itk_3^2}   \\
& = e^{-i t \Delta} \left[\left(e^{it \Delta} u(0) \right)^2 e^{it \Delta} \overline u(0)\right].
\end{align*}
Hence, we find the following approximation of the oscillatory integral in physical space
\begin{equation*}\begin{aligned}
& \int_0^t \mathrm{Osc}(s,\Delta, u(0))\mathrm{d}s  
\\& =t u^2(0) \Big(\varphi_1\left(-2i t\Delta \right)-\varphi_2\left(-2i t\Delta \right)\Big) \overline{u(0)} 
+ t   e^{-i t \Delta} \left[\left(e^{it \Delta} u(0) \right)^2\varphi_2\left(-2i t\Delta \right) e^{it \Delta} \overline u(0)\right] + \mathcal{O}(t^3 \Delta u).
 \end{aligned}
 \end{equation*}
 Plugging the latter into \eqref{Duh4} yields the second-order resonance-based scheme
 \begin{equation}\label{scheme2}
\begin{aligned}
u^{n+1} &= e^{i\tau\Delta} u^n - i  
t e^{i \tau \Delta} \left[(u^n)^2\Big(\varphi_1\left(-2i\tau \Delta \right)-\varphi_2\left(-2i \tau \Delta \right)\Big) \overline{u^n} \right]
\\& -i \tau  \left[\left(e^{i\tau  \Delta} u^n \right)^2\varphi_2\left(-2i \tau \Delta \right) e^{i\tau \Delta} \overline u^n\right] 
  -\frac{\tau^2}{2} \vert u^n\vert^4 u^n 
\end{aligned}
\end{equation}
which was originally introduced in \cite{BS}.

Following the above construction we see that    \eqref{scheme2} introduces a local error at order $ \mathcal{O}\left(t^3  \Delta u \right)$  in polynomials at most of order $\Delta u$, i.e., with loss of two derivatives at most of order two. This is in contrast to classical second-order schemes which require at least the boundedness of four derivatives (see, e.g., \cite{HochOst10,Lubich08}).
\begin{rem}
In order to construct higher-order resonance-based schemes with order $p\geq 3$ we have to  iterate Duhamel's formula \eqref{Duh} up to higher order and discretise in a resonance-based way the appearing iterated integrals up to desired order $p$. The key challenge lies in controlling the  higher order nonlinear frequency interactions  in a structured way. To achieve this decorated trees provide a powerful tool, see  \cite{BS} for a high order framework of resonance-based schemes up to arbitrary order.
\end{rem}

\section{General setting}\label{sec:general}
In Section \ref{sec:idea} we illustrated the main idea of resonances as a computational tool on the example of periodic, cubic NLS equation. In the construction of the schemes we heavily exploited  the periodic boundary conditions and Fourier series expansion which allowed us to explicitly control the resonance structure and nonlinear frequency interactions in the equation. However, it leaves an important question open: {\em Can we  develop resonance-based schemes in a broader setting away from periodic boundary conditions and NLS, i.e., on more general domains, and for a more general class of equations?}

 Indeed we can extend the idea  of resonances as a computational tool to a large class of nonlinear evolution equations and   in many cases  this approach allows  convergence for much rougher data than classical schemes, see for instance \cite{RS,BBS} for the general setting and \cite{BLW,MB,BuyangNS} for various examples such as the Navier--Stokes equation,  parabolic problems with maximum principle, Neumann boundary conditions, etc.  
 
We  illustrate the main idea on the  abstract evolution equation
\begin{equation}\label{ev}
 \partial_t u - \Sigma u = P(u,\overline u) \qquad (t,x) \in \mathbb{R}\times \Omega
\end{equation}
with 
\begin{equation}\label{lapi}
\Sigma \in\{ \Delta , i \Delta\}
\end{equation}
and
$\Omega \subset \mathbb{R}^d$, initial condition
\begin{equation}
\label{init}
u_{/t=0}= u_{0}.
\end{equation}
When $\partial \Omega \neq \emptyset$ we equip the problem with some appropriate homogeneous  boundary conditions. We furthermore assume a polynomial nonlinearity
 \begin{align}\label{poly}
 P(u,\overline u) = u^\ell \overline u^m.
\end{align}

\begin{rem} To present the main idea as clearly as possible we give only formal calculations and focus on the simple case \eqref{ev}.
For the  much more general setting 
\begin{equation}\label{evG}
 \partial_t u - \mathcal{L}u = f(u,\overline u) \qquad (t,x) \in \mathbb{R}\times \Omega
\end{equation}
and rigorous analysis we refer to \cite{RS}, which allows a unified framework  for parabolic,  dispersive as well as mixed equations, covering many examples such as 
\begin{itemize}
\item Nonlinear heat equations
$$ \partial_t u  -  \Delta u =f(u,\overline u),\quad \text{i.e., }\quad \mathcal{L}=\Delta,  \, \mathcal{A}= 0$$
\item Nonlinear Schr\"odinger equations
$$
 i \partial_t u +  \Delta u = \pm \vert u\vert^{2m} u , \,m \in \mathbb{N}\quad \text{i.e., }\quad  \mathcal{L}= i \Delta, \,  \mathcal{A}= -  2 i \Delta , \,     \quad f(u,\overline u) = \pm  i u^{m+1 } \overline u^{m}$$
 \item  Complex Ginzburg Landau equations
\begin{align*} & \partial_t u  -  \alpha  \Delta u = \gamma  u (1 - |u|^2), \,  \alpha, \, \gamma  \in \mathbb{C}, \mbox{Re }\alpha \geq 0, \\  &\text{i.e., }\quad \mathcal{L}=\alpha\Delta, \,  \, \mathcal{A}=  - 2 i \mbox{Im } \alpha \Delta, \, f(u, \overline{u})= \gamma  u (1 - u \overline{u})
\end{align*}
\item Half wave equation
$$ i \partial_t u + \vert \nabla \vert u =  \pm |u|^2 u,\quad \text{i.e., }\quad \mathcal{L}= i \vert \nabla \vert, \,   \quad f(u,\overline u) = \pm u^2 \overline{u}$$
\item  Klein--Gordon and wave-type equations 
$$  \partial_{tt} u -  \Delta u  +  m^2 u =  f(u). $$
\end{itemize}

We could also add potentials, see, e.g., \cite{BBS}.
\end{rem}
To derive the new class of resonance-based schemes for \eqref{ev} we again look at its iteration of Duhamel's formula 
\begin{equation}\label{u1p}
\begin{aligned}
u(t) &= e^{ t \Sigma} u_0 +
\int_0^t e^{ (t-  \xi) \Sigma}  P \left(e^{ \xi\Sigma} u_0,e^{ \xi  \overline\Sigma}\overline u_0\right) d\xi + \int_0^t \int_0^\xi \ldots d\xi_1 d\xi
\end{aligned}
\end{equation}
with (now abstract) leading oscillations
\begin{align}\label{oscG}
\text{Osc}(t,\xi,\Sigma, v,\overline v)
=  e^{ (t-  \xi)\Sigma}  P\left(e^{ \xi\Sigma} v, e^{ \xi \overline\Sigma }\overline v\right) .  
\end{align}
\begin{rem}
Unlike for dispersive equations such as the nonlinear Schr\"odinger equation discussed in Section \ref{sec:idea}, in the general setting of \eqref{ev} we can not split the oscillatory phase $e^{t \Sigma}$ from the oscillations $\text{Osc}$, as $e^{ -  \xi\Sigma}$ might not be well-defined. Hence, to also allow for parabolic type problems (i.e., $\Sigma = \Delta$) we  keep the full term $ e^{ (t-  \xi)\Sigma}$ in \eqref{oscG}
\end{rem}

Thanks to definition \eqref{oscG}  the solution  $u(t)$ defined  in \eqref{u1p} can be expressed as follows
\begin{align}\label{u1}
u(t) = e^{ t \Sigma} u_0 +\int_0^t \text{Osc} (t,\xi,\Sigma, u_0 ,\overline u_0) d\xi+ \mathcal{O}(t^2),
\end{align}
where the remainder $ \mathcal{O}(t^2)$ corresponds to the iterated integral $\int_0^t \int_0^\xi  d\xi_1 d\xi$ and requires no loss of derivatives. In contrast to the periodic case, where we can use Fourier series expansion to explicitly calculate the nonlinear frequency interactions in $\text{Osc} (t,\xi,\Sigma, u_0 ,\overline u_0)$ we here have to take a more general approach to deal with the abstract oscillations \eqref{oscG}: For this purpose we apply the fundamental theorem of calculus and obtain
\begin{equation}\label{osci}
\begin{aligned}
\int_0^\xi \text{Osc} (t,\xi,\Sigma, u_0 ,\overline u_0) d\xi & =  t  \text{Osc}(t, 0,\Sigma, v, \overline{v}) + \int_0^t \int_0^\xi  \partial_s \text{Osc}(t, s, \Sigma, v, \overline{v} ) d s d \xi\\
& 
=  t   e^{ t\Sigma} P(v,\overline v)+ \int_0^t \int_0^\xi  \partial_s\text{Osc} (t, s, \Sigma,v, \overline{v} ) d s d \xi.
\end{aligned}
\end{equation}
Next we calculate that 
\begin{equation}\label{oscil}
\begin{aligned}
& \partial_\xi  \text{Osc}(t,\xi,\Sigma, v,\overline v) \\
&  = e^{(t-\xi)\Sigma} \Bigl[  
\left.  -  \Sigma\left( p \left(e^{\xi\Sigma} v, e^{ \xi \overline\Sigma}  \overline v\right)\right)
+   \ell  \left(e^{\xi\Sigma} v\right)^{\ell-1} \left(e^{ \xi \overline\Sigma}  \overline v\right)^m\cdot\Sigma e^{ \xi\Sigma} v 
+m \left(e^{\xi\Sigma} v\right)^{\ell} \left(e^{ \xi \overline\Sigma}  \overline v\right)^{m-1}\cdot\overline \Sigma e^{ \xi\overline \Sigma} \overline{v} 
 \right]\\
 &  = e^{(t-\xi)\Sigma} \Bigl[  
\left.  
m \left(e^{\xi\Sigma} v\right)^{\ell} \left(e^{ \xi \overline\Sigma}  \overline v\right)^{m-1}\cdot(-\Sigma+\overline \Sigma )e^{ \xi\overline \Sigma} \overline{v} 
 \right] \quad +\quad \text{lower order terms},
\end{aligned}
\end{equation}
where the lower order terms in \eqref{oscil} involve only first-order derivatives, whereas $\text{deg}\,\Sigma =2$ (cf. \eqref{lapi}). 
\begin{rem}
For a general differential operator $\Sigma$ of degree $p$ the lower order terms in \eqref{oscil} will at most be of order $q \leq p-1$.
\end{rem}
\begin{rem}
We   see that for real-valued operators $\Sigma = \overline \Sigma$ the dominant, i.e., leading error term
$$
(-\Sigma+\overline \Sigma )e^{ \xi\overline \Sigma} \overline{v} 
$$ 
in \eqref{oscil} drops. In this case, the classical approximation$$
\int_0^\xi \text{Osc} (t,\xi,\Sigma, u_0 ,\overline u_0) d\xi \approx  t  \text{Osc}(t, 0,\Sigma, v\overline{v}) 
$$
 as taken for instance in  Lie splitting \eqref{Lie}  only involves lower order derivatives. In the general case, $\Sigma$ possibly complex valued, the dominant term
\begin{equation}\label{classAp}
(-\Sigma+\overline \Sigma )e^{ \xi\overline \Sigma} \overline{v} ,
\end{equation}
however, requires the same regularity assumptions as the leading operator $\Sigma$ in the equation. To allow for low regularity approximations (requiring less regularity than the full operator $\Sigma$) in the general setting ($\Sigma$ possibly complex) the classical approximation \eqref{classAp} is no longer sufficient: We have to get rid of the leading error term $(-\Sigma+\overline \Sigma )e^{ \xi\overline \Sigma} \overline{v}$.
\end{rem}
 For this purpose we introduce the resonance-based oscillations
\begin{align}\label{Rosci}
{\text{Rosc}}(t,\xi, r, \Sigma, v,\overline v)
=  e^{ (t-  \xi)\Sigma}  P \left(e^{ \xi\Sigma} v, e^{\xi \Sigma} e^{ r (-\Sigma+ \overline\Sigma) }\overline v\right) 
\end{align}
which satisfy along the diagonal $r=\xi$
$$
{\text{Rosc}}(t,\xi, r=\xi, \Sigma, v,\overline v) =  \text{Osc} (t,\xi,\Sigma, v ,\overline v).
$$
Next, by the fundamental theorem of calculus we compute that
\begin{align}\label{Rosc:exp}
 \text{Osc} (t,\xi,\Sigma, v ,\overline v) = {\text{Rosc}}(t,0, \xi, \Sigma, v,\overline v) 
 + \int_0^\xi \partial_s{\text{Rosc}}(t,s, \xi, \Sigma, v,\overline v) ds,
\end{align}
where
\begin{equation}\label{oscilR}
\begin{aligned}
& \partial_s{\text{Rosc}}(t,s, \xi, \Sigma, v,\overline v) \\
&  = e^{(t-s)\Sigma} \Bigl[  
\left.  -  \Sigma\left( P \left(e^{s\Sigma} v, e^{s\Sigma} e^{ \xi (-\Sigma+ \overline\Sigma) }  \overline v\right)\right)
+   \ell  \left(e^{s\Sigma} v\right)^{\ell-1} \left(e^{s\Sigma} e^{ \xi (-\Sigma+ \overline\Sigma) }    \overline v\right)^m\cdot\Sigma e^{ s\Sigma} v \right. 
\\&\left.\qquad + \,m \left(e^{s\Sigma} v\right)^{\ell} \left(e^{s\Sigma} e^{ \xi (-\Sigma+ \overline\Sigma) }    \overline v\right)^{m-1}\cdot \Sigma e^{s\Sigma} e^{\xi (-\Sigma+ \overline\Sigma) }   \overline{v} 
 \right]\\
 &  =\quad  \text{lower order terms only involving differential operators of degree $q \leq \text{deg}(\Sigma)-1$}.
\end{aligned}
\end{equation}
This is the main advantage of the new resonance-based approach \eqref{Rosci}: Instead of introducing the full error term \eqref{oscil}, the improved remainder \eqref{oscilR} only involves lower order differential operators.\smallskip

Next we plug the resonance-based expansion \eqref{Rosc:exp} together with the observation that (cf. \eqref{Rosci})
$$
{\text{Rosc}}(t,0, \xi, \Sigma, v,\overline v)
= e^{ t \Sigma}  P \left(  v,   e^{\xi (-\Sigma+ \overline\Sigma) }\overline v\right) 
$$
into the iteration of Duhamel's formula \eqref{u1} and obtain
\begin{align*}
u(t) = e^{ t \Sigma} u_0 +\int_0^t e^{ t \Sigma}  P \left(  u_0,   e^{\xi (-\Sigma+ \overline\Sigma) }\overline u_0\right) 
   d\xi+ t^2 \mathcal{R}_ { \text{deg}(\Sigma)-1}( u),
\end{align*}
where the remainder $\mathcal{R}_ { \text{deg}(\Sigma)-1}( u)$ involves only lower order differential operators of degree $ \text{deg}(\Sigma)-1$. Exploiting the polynomial structure of  nonlinearity $P$ (see \eqref{poly}) we furthermore obtain  that
\begin{align}\label{u1Rosc}
u(t) = e^{ t \Sigma} u_0 + t e^{ t \Sigma}  P \left(  u_0,   \varphi_1(t (-\Sigma+ \overline\Sigma) )\overline u_0\right) +  t^2 \mathcal{R}_ { \text{deg}(\Sigma)-1}(  u).
\end{align}
This builds  the motivation for the general first-order resonance-based scheme
\begin{align*}
u^{n+1} = e^{ \tau \Sigma} u^n + \tau e^{ \tau \Sigma}  P \left(  u^n,   \varphi_1(\tau (-\Sigma+ \overline\Sigma) )\overline u^n\right)
 \end{align*}
 which allows convergence in a general setting, including parabolic, dispersive and hyperbolic problems, for rougher data than classical schemes. For a precise convergence analysis we refer to \cite[Theorem 2]{RS}.
 
 \begin{rem}
 Similarly to the periodic setting of Section \ref{subsec:rc} we can also achieve higher order approximations in the general setting, by considering higher order Duhamel iterations in \eqref{u1p} and their corresponding resonance-based discretisation. For a high-order analysis in the general setting via decorated tree series analysis we refer to \cite{BBS}.
 \end{rem}

\begin{rem}Another natural question is to see in how far resonances can be used when randomness comes into play. We refer to \cite{BBS23} and \cite{Ystoch} for recent progress in this direction.
\end{rem}

\section{Error analysis: or how far can we go?}\label{sec:error}
In this section we want to know how far we can actually push down the Sobolev index $\sigma>0$ when solving dispersive equations such as cubic NLS (cf. \eqref{iniSig}) with the fundamental question: can we prove convergence  at the level of well-posedness of the equation? More precisely, if we know that the solution of the equation exists  in some Sobolev space $H^{\gamma_1}$, $\gamma_1>0$ (globally in time or up to a certain time $T>0$) can we then also approximate the solution at this level of regularity, i.e., establish an error estimate
$$
\Vert u(t) - u^n\Vert_{H^{\gamma_1}} \leq \tau^q  C_{T},
$$
where $C_{T}$ depends only on the $H^{\gamma_2}$ norm of the solution $u$ on $[0, T]$ and  the exponent $q$ naturally will depend on $\gamma_1$ and $\gamma_2$ and we wish to choose $\gamma_2 >\gamma_{1}$ as close as possible to $\gamma_1$.

This question becomes in particular tough to answer when we indeed want to go below the critical regularity
\begin{equation}\label{crit}
0 <\gamma_2\leq d/2
\end{equation} where classical techniques in error analysis based on standard product estimates break down.

In order to explain the ideas, we focus on  the cubic, periodic NLS \eqref{cubNLS} in dimension one, though the techniques we shall introduce can be used
to perform error estimates at low regularity of time discretizations for  a large class of dispersive equations with periodic boundary conditions
(see for example \cite{ORS22}, \cite{ORSsplit} and  \cite{2DNLS} for higher dimensions,  \cite{RSsplit} for the case of the KdV equation).
These ideas are usefull to analyze not only the resonance based schemes presented in the previous section but also more
classical exponential integrators or splitting schemes.
To understand the difficulties for the error analysis at low regularity,  we firsts  have to recall 
the main difficulties that show up when one tries to prove  local well posedness in $H^s$, $s<1/2$ at the continuous level.
In the whole space $x\in \mathbb{R}$ the analysis relies  on  Strichartz estimates which are valid because of the dispersive decay properties (see for example the book \cite{Tao06}): one can prove that for the linear Schrodinger equation
\begin{equation}
\label{linSchro} i\partial_{t} u + \partial_{xx} u = F, \, x \in \mathbb{R}, \quad u_{/t=0}= u_{0}
\end{equation}
one gets the estimate
$$ \|u \|_{L^p_{t}(\mathbb{R}, L^q_{x}(\mathbb{R}))} \lesssim_{p,q} \|u_{0}\|_{L^2} + \|F\|_{L^{p_{1}'}_{t}(\mathbb{R}, L^{q_{1}'}_{x}(\mathbb{R}))}$$
 where $(p,q)$ and $(p_{1}, q_{1})$ are admissible pairs and $p_{1}'$, $q_{1}'$ are such that $p_{1}^{-1}+ (p_{1}')^{-1}=1$, 
 $q_{1}^{-1}+ (q_{1}')^{-1}=1$. The pair $(p,q)$ is admissible (in dimension one)  if $p \geq 2$ and
  $ 2/p + 1/q= 1/2$.  In particular one gets that for an $L^2$ initial data, the solution $u$ of \eqref{linSchro} is in $L^4_{t} L^\infty_{x}$.
   This strong gain of integrability allows to prove local existence for rough data by a fixed point argument.
   This estimate strongly relies on the dispersive estimate: there exists $C>0$ such that for every solution of \eqref{linSchro}
   with $F=0$, one has
   $$ \|u(t)\|_{L^\infty(\mathbb{R})} \leq { C \over t^{1\over 2}} \|u_{0}\|_{L^1(\mathbb{R})}, \quad \forall t>0.$$
{Up to our knowledge in numerical analysis discrete type Strichartz estimates were first applied to cubic Schr\"odinger equations with solutions in $H^2(\mathbb{R}^d)$, $d \leq 3$ in  \cite{IZ06,IZ09,Ignat11}}.
    
This estimate is not true for periodic boundary conditions and thus the Strichartz estimates do not hold in such generality.
Nevertheless, for periodic boundary conditions, there is still  a gain of integrability which holds locally in time.
 For example, for  the linear Schr\"odinger equation 
 $$ i \partial_{t} u + \partial_{x} u = 0, \, x \in \mathbb{T},  \quad u_{/t=0}=u_{0},$$
 we have
 \begin{equation}
 \label{strichtorus} \| u \|_{L^4(\mathbb{T}^2)} \lesssim \|u_{0}\|_{L^2}.
 \end{equation}
 This classical estimate in harmonic analysis goes back to Zygmund. A useful way to encode this type of properties to analyze
 nonlinear problems with periodic boundary conditions is to use the Fourier restriction spaces introduced by Bourgain \cite{Bour93a}.
 For a function $u(t,x)$ on $\mathbb{R}\times\mathbb{T}$, $\tilde{u}(\sigma,k)$ stands for its time-space Fourier transform, i.e.
$$
\tilde{u}(\sigma,k)=\int_{\mathbb{R}\times\mathbb{T}}u(t,x)e^{-i\sigma t-i kx}dx dt.
$$
The inverse transform is given by
$$
u(t,x)=\sum\limits_{k\in\mathbb{Z}}\hat{u}_k(t)e^{i kx}
$$
with the Fourier coefficients $\hat{u}_k(t)=\frac{1}{2\pi}\int_{\mathbb{R}}\tilde{u}(\sigma,k)e^{i\sigma t}d\sigma$.

This way, we can define the Bourgain space $X^{s, b}= X^{s,b}(\mathbb{R}\times\mathbb{T})$ consisting of tempered distribution with finite norm
$$
\Vert u\Vert_{X^{s,b}}=\Vert\langle k\rangle^s\langle\sigma + k^2\rangle^b\tilde{u}(\sigma,k)\Vert_{L^2 l^2}.
$$

Some well-known basic properties of these spaces are the following (we refer for example to \cite{Bour93a}, and the books \cite{Linares-Ponce}, \cite{Tao06}):
\begin{lemma}
\label{lemmabourgainfacile}
For $\eta\in \mathcal{C}^\infty_{c}(\mathbb{R})$, we have that
\begin{eqnarray}
\label{bourg1an}
&  \| \eta (t) \e^{i t \partial_x^2 } f \|_{X^{s,b}} \lesssim_{\eta, b} \|f\|_{H^s}, \quad s \in \mathbb{R},\, b\in \mathbb R,\, f \in H^s(\mathbb{T}), \\
  \label{bourg2an}
&   \| \eta (t) u  \|_{X^{s,b}} \lesssim_{\eta, b} \|u\|_{X^{s,b}}, \quad s \in \mathbb{R}, \, b\in\mathbb R, \\
   \label{bourg3an}
&  \| \eta({t \over T}) u \|_{X^{s,b'}} \lesssim_{\eta, b, b'} T^{b-b'} \|u\|_{X^{s,b}}, \quad s \in \mathbb{R},  -{1 \over 2} <b' \leq b <{ 1 \over 2}, \, 0< T \leq 1, \\
& \left\| \eta(t) \int_{-\infty}^t \e^{i (t-s) \partial_{x}^2} F(s) \dd s \right\|_{X^{s,b}} \lesssim_{\eta, b} \|F\|_{X^{s, b-1}}, \quad s \in \mathbb{R}, \,
 b>{1 \over 2}, \\
 &  \label{bourg4an}  \|u\|_{L^\infty(\mathbb{R}, H^s)} \lesssim_{b} \|u\|_{X^{s, b}}, \quad b>1/2, \, s \in \mathbb{R}.
\end{eqnarray}
\end{lemma}
We actually have  the continuous embedding $X^{s, b} \subset \mathcal{C}(\mathbb{R}, H^s)$ for $b>1/2$ which is a consequence
of the Sobolev embedding in the time variable.
 
 The crucial estimate for the analysis of the cubic equation  is the following appropriate generalization of \eqref{strichtorus}
  in the framework of Bourgain spaces:
  \begin{lemma}
\label{lemmabourgaindur}
There exists a constant $C>0$ such that for every $u \in X^{0, {3\over 8}}$, we have the estimate
$$
\| u \|_{L^4( \mathbb{R} \times \mathbb{T})} \leq C \|u \|_{X^{0, {3 \over 8}}}.
$$
\end{lemma}
Again, we refer  to \cite{Tao06} Proposition 2.13  for its proof. Note that, by duality, we also obtain that
$$
\| u \|_{X^{0, {-{3 \over 8}}}} \lesssim \| u \|_{L^{4\over 3}( \mathbb{R} \times \mathbb{T})}.
$$
By combining the two estimates with H\"{o}lder, this further implies that
\begin{equation}
\label{prodbourgain}
\| u v w  \|_{X^{0, {-{3 \over 8}}}} \lesssim \|u\|_{X^{0, {3 \over 8}}}  \|v\|_{X^{0, {3 \over 8}}}   \|w\|_{X^{0, {3 \over 8}}}.
\end{equation}

The well-posedness for \eqref{cubNLS} at low regularity  is proven by looking for a fixed point of the functional
\begin{equation}\label{fixed}
F(v)(t)=   \eta(t)\e^{i t \partial_{x}^2} u_{0} -  i \eta(t)\int_{0}^t \e^{i(t- s)\partial_{x}^2 }
\left( \eta\left({s \over \delta}\right) |v(s)|^2 v(s)\right)\dd s,
\end{equation}
where $\eta \in [0, 1]$ is a  smooth compactly supported function which is equal to $1$ on $[-1, 1]$ and supported in $[-2, 2]$. For $|t| \leq \delta \leq 1/2 $, a fixed point of the above functional gives a solution of the original Cauchy problem, denoted by $u$.

By choosing $b \in (1/2, 5/8)$, we can use the estimates of the two previous Lemma to get that 
$$
\|F(v) \|_{X^{0, b}}  \leq C \|u_{0}\|_{L^2} + C \delta^{ \epsilon_{0}}  \|v\|_{X^{0,b}}^3
$$
and if  $\|v_{1}\|_{X^{0, b}} \leq R, \, \|v_{2} \|_{X^{0, b}}\leq R$, that
$$
\|F(v_{1}) - F(v_{2}) \|_{X^{0, b}} \leq 4 C \delta^{\epsilon_{0}} R^2 \|v_{1} - v_{2} \|_{X^{0, b}}
$$
where
$\epsilon_{0}= 5/8 - b >0$.

One can thus use Banach fixed point Theorem in a suitable ball of $X^{0, b}$.
 This gives short time  well posedness for   $ L^2$ initial data. One then get a global solution by using the conservation of the $L^2$ norm.
 This type of estimates also allows to propagate higher regularity globally in time.
 \newline
 
 One can then use these ideas to perform error analysis for a large class of  time discretizations by defining time discrete
 Bourgain spaces.

For a sequence $(u_{n})_{n \in \mathbb{Z}}$, we define its Fourier transform as
$$
\mathcal{F}_{\tau}(u_{n}) (\sigma)= \tau \sum_{m \in \mathbb{Z}} u_{m} \,\e^{ i m \tau \sigma}.
$$
This defines a periodic function on $ [-\pi/\tau, \pi/\tau]$ and we have the inverse Fourier transform formula
$$
u_m =  \frac1{2\pi} \int_{-{\pi \over \tau}}^{\pi \over \tau}  \mathcal{F}_{\tau}(u_{n})(\sigma) \,\e^{- i m \tau \sigma}\dd \sigma.
$$
With these definitions the Parseval identity reads
$$
\|u_{n}\|_{l^2_{\tau}}= \| \mathcal{F}_{\tau} (u_n) \|_{L^2(-\pi/\tau, \pi/\tau)},
$$
where the norms are defined by
$$
\|u_{n}\|_{l^2_{\tau}}^2 = \tau \sum_{n \in \mathbb{Z}} |u_{n}|^2, \quad
\| \mathcal{F}_{\tau} (u_n) \|_{L^2(-\pi/\tau, \pi/\tau)}^2 = {1 \over 2 \pi} \int_{-{\pi \over \tau}}^{\pi \over \tau} | \mathcal{F}_{\tau} (u_n) (\sigma) |^2
\dd \sigma.
$$
We write $L^2$ instead of $L^2 (-\pi/\tau, \pi/\tau)$ for short. We stress the fact that this is not the standard way of normalizing the Fourier series.

We then define in a natural way Sobolev spaces $H^b_{\tau}$ of sequences $(u_{n})_{n\in\mathbb{Z}}$ by \\
$$
\| u_{n}\|_{H^b_{\tau}} =  \left\| \langle d_{\tau}( \sigma)\rangle ^b \mathcal F_\tau(u_n) \right \|_{L^2},
$$
with $d_{\tau}(\sigma)=\frac{\e^{i \tau \sigma} - 1}\tau$ so that we have equivalent norms
$$
\| u_{n}\|_{H^b_{\tau}} =\|\langle D_{\tau}\rangle^b u_{n}\|_{l^2_{\tau}},
$$
where the operator $D_{\tau}$ is defined by $ (D_{\tau}(u_{n}))_{n}=\left( \frac{u_{n-1} - u_{n}}\tau\right)_{n}$  since by definition of the Fourier transform
$$
\mathcal{F}_{\tau}( D_{\tau} u_{n}) (\sigma) =d_{\tau}(\sigma) \mathcal{F}_{\tau}(u_{n}) (\sigma).
$$
Note that $d_{\tau}$ is $2\pi/\tau$ periodic and that uniformly in $\tau$, we have $|d_{\tau}(\sigma)| \sim | \sigma |$ for $|\tau \sigma | \leq \pi$.

For sequences of functions $(u_{n}(x))_{n \in \mathbb{Z}},$ we define the Fourier transform $\widetilde{u_{n}}(\sigma, k)$ by
$$
\mathcal F_{\tau,x}(u_n)(\sigma,k) =\widetilde{u_{n}} (\sigma, k)= \tau \sum_{m \in \mathbb{Z}} \widehat{u_{m}}(k) \,\e^{i m \tau \sigma}, \quad \widehat{u_{m}}(k)= {1 \over 2\pi} \int_{-\pi}^\pi u_{m}(x) \,\e^{-i k x}\dd x.
$$
Parseval's identity then reads
\begin{equation}\label{parseval}
\| \widetilde{u_{n}}\|_{L^2l^2}= \|u_{n}\|_{l^2_{\tau}L^2},
\end{equation}
where
$$
\| \widetilde{u_{n}}\|_{L^2l^2}^2 = \int_{-{\pi \over \tau}}^{\pi\over \tau} \sum_{k \in \mathbb{Z}}
|\widetilde{u_{n}}(\sigma, k)|^2 \dd \sigma, \quad
\|u_{n}\|_{l^2_{\tau}L^2}^2 = \tau \sum_{m \in \mathbb{Z}} \int_{-\pi}^\pi  |u_{m}(x)|^2 \dd x.
$$
We then finally  define the discrete Bourgain spaces $X^{s,b}_\tau$ for $s\in \mathbb{R}$, $b\in\mathbb R$, $\tau>0$ by
\begin{equation}\label{bourgdef1}
\| u_n \|_{X^{s,b}_{\tau}}= \| \e^{-i n \tau \partial_{x}^2 } u_{n}\|_{H^b_{\tau} H^s }=  \|
\langle D_{\tau}\rangle^b \langle \partial_{x} \rangle^s (\e^{-i n \tau \partial_{x}^2 } u_{n})\|_{l^2_{\tau} L^2 }
\end{equation}
where $\langle \cdot \rangle= (1 + |\cdot|^2)^{1 \over 2}$.

As in the continuous case, we can relate this  definition to a Fourier restriction norm adapted to the dispersion relation
of the underlying scheme: 

\begin{lemma}\label{lemequiv}
With the above definition, we have that
\begin{equation}\label{norm2}
\| u_n \|_{X^{s,b}_{\tau}} \sim  \left\| \langle k \rangle^s \langle  d_{\tau}(\sigma -k^2)  \rangle^b \widetilde{u_n}(\sigma, k)  \right\|_{L^2l^2}.
\end{equation}
\end{lemma}
Note that the weight $d_{\tau}(\sigma -k^2)$ vanishes if $\tau(\sigma - k^2)= 2 m \pi$ for $m \in \mathbb{Z}$. For a localized function such that $k$ is constrained to $|k| \lesssim  \tau^{-{1 \over 2}}$ this will behave like in the continuous case with only a cancellation when $\sigma = k^2$. For larger frequencies, however, there are additional cancellations that will create some loss in the product estimates at the discrete level.

%


With this definition, the counterpart of Lemma \ref{lemmabourgainfacile} holds  at the discrete level.

\begin{lemma}\label{bourgainfaciled}
For $\eta \in \mathcal{C}^\infty_{c}(\mathbb{R})$ and $\tau\in(0,1]$, we have that
\begin{align}
\label{bourg1} &\| \eta(n \tau)  \e^{in \tau \partial_{x}^2} f\|_{X^{s,b}_{\tau}} \lesssim_{\eta, b} \|f\|_{H^s}, \quad s \in \mathbb{R}, \, b \in \mathbb{R}, \, f \in H^s, \\
\label{bourg2} &\| \eta(n \tau)  u_{n}\|_{X^{s,b}_{\tau}} \lesssim_{\eta, b} \|u_{n}\|_{X^{s,b}_{\tau}}, \quad s \in \mathbb{R}, \, b \in \mathbb{R} , \, u_{n} \in X^{s,b}_{\tau},\\
\label{bourg3} &\left\| \eta\left(\frac{n\tau}T \right) u_{n} \right\|_{X^{s,b'}_{\tau}} \lesssim_{\eta, b, b'} T^{b-b'} \|u_{n}\|_{X^{s,b}_{\tau}}, \quad s \in \mathbb{R},  -{1 \over 2} <b' \leq b <{ 1 \over 2},\, 0< T = N \tau  \leq 1, \, N \geq 1.
\end{align}
In addition, for
$$
U_{n}(x)= \eta(n \tau) \tau \sum_{m=0}^n  \e^{i ( n-m ) \tau \partial_{x}^2}  u_{m}(x),
$$
we have
\begin{equation}
\label{bourg4}\|U_{n}\|_{X^{s,b}_{\tau}} \lesssim_{\eta, b} \|u_{n}\|_{X^{s, b-1}_{\tau}}, \quad s \in \mathbb{R}, \, b>1/2.
\end{equation}
We stress that all given estimates are uniform in $\tau$.
\end{lemma}

Finally, we have to study  the discrete counterpart of Lemma \ref{lemmabourgaindur} which is crucial for the analysis of nonlinear problems at low regularity.

In the discrete setting, for a sequence $(u_{n}) \in l^p(\mathbb{Z}, X)$, with $X$ normed space  we use the norm
\beq
\label{normelptau} \|u_{n}\|_{l^p_{\tau} (X)}=\left( \tau \sum_{n \in \mathbb{Z}} \|u_{n}\|_{X}^p \right)^{1 \over p}.
\eeq
We shall use  the Fourier multiplier 
$$ \widehat{\Pi_{K} f} = \mathrm{1}_{|k| \leq K} \hat{f}, \quad f \in L^2(\mathbb{T})$$
which projects on frequencies smaller than $K$.

\begin{lemma}\label{prodd}
For $K \geq \tau^{-{ 1 \over 2}}$,  we have
\begin{equation}\label{prodd1}
\|\Pi_{K}u_{n} \|_{l^4_{\tau}L^4} \lesssim  (K \tau^{ 1 \over 2})^{1 \over 2} \|u_{n}\|_{X^{0, {3\over 8}}_{\tau}}.
\end{equation}
\end{lemma}

Note that choosing $K= \tau^{-{1 \over 2}}$, we get the same estimate as in the continuous case. This can be interpreted
as a type of CFL condition. If we impose only a weaker condition then the properties of the wave interactions at the discrete level
are different from the ones at the continuous level, in some sense, the interaction picture of the continuous case  is reproduced $K \tau^{- {1 \over 2}}$ times in the discrete case and we start loosing uniformity. By using a Littlewood-Paley decomposition, 
this can be translated into a loss of derivative in the embedding. We can get the following
\begin{cor}
For $K= \tau^{- {\alpha \over 2}}$, $\alpha \geq 1$, we have 
$$ \|\Pi_{K}u_{n} \|_{l^4_{\tau}L^4} \lesssim   \|u_{n}\|_{X^{{1 \over 2}(1- {1 \over \alpha}), {3\over 8}}_{\tau}}$$
\end{cor}
Note that for $\alpha <2$ this estimate is still an improvement compared to the Sobolev embedding $H^{1 \over 4}(\mathbb{T}) \subset L^4(\mathbb{T}).$
\newline.

To illustrate how this framework can be used to perform energy estimate,  we shall study error estimates for 
the following filtered splitting scheme
\begin{equation}\label{scheme}
\begin{aligned}
u^{n+1}  &  = \e^{i \tau \partial_x^2}  \Pi_{\taul}\left( \e^{-i \tau \vert \Pi_{\taul} u^n \vert^2} \Pi_{\taul} u^n \right),
\qquad u^0 & =\Pi_{\taul}  u(0)
\end{aligned}
\end{equation}
where we have set $\Pi_{\tau}= \Pi_{K}$ with $K= \tau^{-{1 \over 2}}$ for short. In view of the above analysis, 
the presence of the filter is crucial to have at our disposal optimal low regularity product estimates.
The filtered Lie splitting \eqref{scheme} can be seen as a classical Lie splitting discretisation of the projected equation
\begin{equation}
\begin{aligned}\label{nlsP}
i \partial_t  \uta & = - \partial_{x}^2 \uta  +  \Pi_{\taul}\left( \vert  \Pi_{\taul}   \uta \vert^2  \Pi_{\taul}  \uta\right), \qquad \uta(0) & =  \Pi_{\taul}  u(0).
\end{aligned}
\end{equation}
For an initial data $u_{0} \in H^{s_{0}}$, and for every $T >0$
by using the PDE estimates, one can first show that
 uniformly for $\tau \in (0, 1]$, we have for some $C_{T}>0$,
$$
\| u - u_{\tau} \|_{L^\infty((0,T); L^2)} \lesssim \| u - u_{\tau} \|_{X^{0, b}(T)} \leq  C_{T}  \tau^{s_{0} \over 2}
$$
where $X^{0, b}(T)$ is the local version of  the Bourgain space.
This yields that in order to get an error estimate, it suffices to study $e_{n} = u_{n} -u_{\tau} (t_{n})$. 
As above, we then study an extended equation for the error under the form
\begin{equation}
\label{eqen}
e^{n}=- i  \tau \eta(t_{n})   \sum_{k=0}^{n-1}  \e^{i (n-k)  \tau \partial_x^2} \eta\left( {t_{k} \over T_{1}}\right)
\bigl( \Phi_{\mathcal{N}}^\tau(\vta(t_k)) - \Phi_{\mathcal{N}}^\tau(u_{\tau}(t_{k})-e^k)\bigr) + \mathcal{R}_{n},
\end{equation}
\begin{equation}
\label{Rn}
\mathcal{R}_{n}=- i \eta(t_{n})  \sum_{k=0}^{n-1}  \e^{i (n-k)  \tau \partial_x^2} \eta(t_{k}) \mathcal{E}_{\text{loc}}(t_k,\tau,\vta)
\end{equation}
where $ \mathcal{E}_{\text{loc}}(t_k,\tau,\vta)$ is the local error for the time discretization of $u_{\tau}$ and  we have set 
$$\Phi_{\mathcal{N}}^\tau(w) =- \Pi_{\taul} \left( \frac{\e^{-i \tau \vert \Pi_{\taul} w \vert^2}-1}{i  \tau} \Pi_{\taul} w\right).$$
Note that for  $0 \leq n \leq N_{1}$, where $N_{1}=\lfloor {T_{1}\over \tau}\rfloor$ with $T_{1} \leq \min(1, T)$, 
the solution of \eqref{eqen} is indeed the error we want to estimate.

The analysis of the local error yields
 $$\|\mathcal{R}_{n}\|_{X^{0, b}_{\tau}}\leq C_{T}\tau^{s_{0} \over 2}.$$
 Moreover, by using Lemma \ref{bourgainfaciled} and Lemma \ref{prodd}, we can obtain an estimate under the form  
$$ \|e^n \|_{X^{0, b}_{\tau}} \leq   C_T T_{1}^{\varepsilon_{0}}
 \left(  \|e^n \|_{X^{0, {3 \over 8}}_{\tau}} + \|e^n \|_{X^{0, {3 \over 8}}_{\tau}}^2 + \|e^n \|_{X^{0, {3 \over 8}}_{\tau}}^3\right) +
   C_{T}\tau^{ s_{0}\over 2}$$
 for some $\varepsilon_{0}>0$  for the solution of \eqref{eqen}.
 
 This yields the following result
 \begin{theorem}\label{maintheo}
For every $T >0$ and $u_{0} \in H^{s_{0}}$,  $s_{0}>0$, let $u \in \mathcal{C}([0, T], H^{s_{0}})$ be the exact solution of \eqref{cubNLS} with initial datum $u_{0}$ and denote by $\un$ the sequence defined by the scheme~\eqref{scheme}. Then, we have  the following error estimate:
there exists $\tau_{0}>0$ and $C_{T}>0$ such that for every step size $\tau \in (0, \tau_{0}]$
\beq
\label{final1} \|\un- u(t_{n})\|_{L^2} \leq C_{T} \tau^{s_{0}\over 2}, \quad 0 \leq n\tau \leq T.
\eeq
\end{theorem}
We refer to \cite{ORSsplit}, \cite{2DNLS} for the full proof.

Here we only discuss the temporal error; for fully discrete error estimates, see \cite{2DNLSd}.
\section{Open questions}\label{sec:open}
When it comes to physically motivated equations a key challenge in numerical integration lies in structure preservation: If the equation has a certain balance law (e.g., energy or mass conservation), we would like to reproduce this structure as far as possible also at the numerical (discrete) level. For ODEs an impressive theory could be meanwhile established with a {revolutionary step}   set by the theory of  {geometric numerical integration} (cf. Hairer et al. \cite{HLW}, Blanes $\&$ Casas \cite{BlanesGNI}, Engquist et al. \cite{EFHI09},    Leimkuhler  $\&$ Reich  \cite{LR04}, Sanz-Serna $\&$ Calvo \cite{SanBook}).  
{\color{black} For PDEs, on the other hand, a unified   theory is far out of reach and very little is known so far}, see, e.g., the groundbreaking and pioneering works on the analysis of splitting and trigonometric methods for nonlinear Schr\"odinger and wave equations  over long times \cite{CHL08a,GLF,GL08b,GL10,Faou12}, scattering \cite{Scat}, long-time error estimates for small, smooth initial data \cite{Bao,Long1,Long2} and energy and mass preserving Crank--Nicolson approximations for the nonlinear Schr\"odinger equation   \cite{Akrivis,Besse21,PetCN,SanNLS}.

It is a natural, yet in large parts widely open question to ask in how far resonances can be used to design structure preserving schemes.  First we have to note that classical resonance-based schemes are in general not structure preserving,  they  are even  not   symmetric. 
 \begin{ex}[Loss of symmetry in classical resonance-based schemes] \vskip-0.2cm\label{exNotS}
 We say that a one-step method $u^{n+1} = \Phi_\tau (u^n)$  is symmetric, if it holds  that 
$$
\Phi_\tau = \Phi_{-\tau}^{-1}.
$$
Let us for example consider   the second order resonance-based scheme  \eqref{scheme2} for cubic Schr\"odinger equation. While the scheme~\eqref{scheme2}  allows us to reduce the regularity assumptions imposed by classical schemes  it destroys symmetry. This can be easily seen by the following observation: 
The adjoint method $\Phi_{-\tau}^{-1}$ is implicit, while the scheme~\eqref{scheme2} itself is explicit. The scheme \eqref{scheme2} cannot therefore  be symmetric.
\end{ex}

This is also what we observe in  simulations: energy and mass are not well preserved over long time scales. Recently, symmetric resonance-based schemes could be introduced in \cite{YYMS}.  Symplectic first-order resonance-based schemes for periodic, one-dimensional NLS and KdV could be found in \cite{MS23}. A general class of symplectic resonance-based schemes  (beyond one dimension and two special cases) which inherit the symplectic structure of the underlying PDE  on the other hand remain open.  The main difficulty lies in the open question how to find a {\em symplectic} resonance-based discretisation of the leading oscillation, which are in case of cubic NLS posed on $\mathbb{T}^d$ for instance driven by
 \begin{equation*}\label{Res:NLS}
\mathrm{R}( {\bf k_1}, {\bf k_2}, {\bf k_3}) =   2 {\bf k_2}\cdot - 2 ( {\bf k_1}+ {\bf k_3}) \cdot  {\bf k_2} + 2  {\bf k_2} \cdot  {\bf k_3}, \quad  {\bf k_1},  {\bf k_2},  {\bf k_3} \in \mathbb{Z}^d.
 \end{equation*}

\bibliographystyle{alpha}

\end{document}